\documentclass[12pt]{article}

\title{A Natural Basis of States for the Noncommutative Sphere and its
Moyal bracket}
\author{{\bf J. Gratus}%
\\
Laboratoire de Gravitation et Cosmologie Relativistes%
\thanks{\it Laboratoire associ\'e au CNRS {\rm URA 769}}
\\
Tour 22/12 4eme etage, Boite Courrier142, 4pl Jussieu. F75252
Paris
\\
email: gratus@ccr.jussieu.fr
}

\date{March 17, 1997} 

\usepackage{amsfonts}
\usepackage{ifthen}

\def\Cmpx{{\mathbb{C}}}
\def\Real{{\mathbb{R}}}
\def\Intg{{\mathbb{Z}}}

\def\sign{{\hbox{\rm sign}}}
\def\degreerm{{\mathrm {degree}}}
\def\ad{\hbox{\rm ad}}
\def\tr{\hbox{\rm tr}}
\def\ftr{\hbox{\rm T}\!{\rm r}}

\def\Pexpr{{\cal P}}
\def\Jexpr{{\cal J}}
\def\Sexpr{{\cal U}}
\def\calS{{\cal S}}
\def\Czz{{C_{\scriptscriptstyle{00}}}}

\def\cnj{\overline}
\def\fs#1{{\mbox{\boldmath $#1$}}}

\def\tfrac#1#2{{\textstyle{\frac{#1}{#2}}}}
\def\nfrac#1#2{{\raisebox{.5ex}{$#1$}\!/\!\raisebox{-.5ex}{$#2$}}}
\def\scrhalf{{\raisebox{.3ex}
{$\scriptstyle 1$}\!/\!\raisebox{-.3ex}{$\scriptstyle 2$}}}

\def\lgap{\hbox{\vrule height 1.6em width 0em}}

\def\Ssym(#1,#2,#3){{S(#1,#2,#3)}} 

\def\vvec(#1,#2){{|#1,#2\rangle}} 
\def\lvec(#1,#2){{\langle#1,#2|}}

\newtheorem{theorem}{Theorem}
\newtheorem{lemma}[theorem]{Lemma}
\newtheorem{corol}[theorem]{Corollary}

\def\ppmatrix#1{\pmatrix{ #1 }}

{%
\refstepcounter{theorem}%
{\vskip 1em\noindent \bf Example\ \thetheorem :}}%
{}

\newenvironment{proof}[1]%
{\vskip 1em\noindent \nopagebreak{\bf 
\ifthenelse{\equal{#1}{}}{Proof:}{Proof #1:}\par}}%
{\hfill\vrule height 0.7 em width 0.7 em\vskip 1em}

\begin{document}

\parskip 0em
\parindent 3em
\abovedisplayshortskip 0pt
\belowdisplayshortskip 3pt
\abovedisplayskip 3pt
\belowdisplayskip 3pt

\maketitle

\begin{abstract}
An infinite dimensional algebra which is a non-decomposable reducible
representation of $su(2)$ is given. This algebra is defined with
respect to two real parameters. If one of these parameters is zero the
algebra is the commutative algebra of functions on the sphere,
otherwise it is a noncommutative analogue. This is an extension of the
algebra normally refered to as the (Berezin) quantum sphere or
``fuzzy'' sphere.  A natural indefinite ``inner'' product and a basis
of the algebra orthogonal with respect to it are given. The basis
elements are homogenious polynomials, eigenvectors of a Laplacian, and
related to the Hahn polynomials.  It is shown that these elements tend
to the spherical harmonics for the sphere. A Moyal bracket is
constructed and shown to be the standard Moyal bracket for the sphere.
\end{abstract}

\newpage



\section{ Introduction }
\label{ch_intro}

The noncommutative or ``fuzzy'' sphere has been considered by several
authors in different contexts. It is an example for a general
quantisation procedure \cite{Berezin74,Berezin75a,CaheGutt90}. 

It is also an example often used in noncommutative geometry
\cite{Madore2,Madore_book,klimcik1,Watamura1} (see also references
within), and in the theory of membranes \cite{DeWit1} which has
application to supersymmetry. It is studied in relation to coherent
states, \cite{Grosse1}, \cite{Perelomov1} and as a reduction of the
symplectic algebra on $\Real^6$ \cite{Molin1} \cite{Molin2}.

Normally the approximation for the algbra of functions on a sphere is
in terms of matrices, where the functions on a sphere appear only in
the limit as the size of the matrix tends to infinity.  In this
article, however, we examine a two parameter algebra of polynomials
$\Pexpr(\kappa,R)$ with $\kappa,R\in\Real$.  For different values of
$\kappa$ and $R$ we obtain:
\begin{trivlist}
\item\qquad $\bullet$
The commutative algebra of finite
sums of harmonics on the sphere (when $\kappa=0$).
In this case $R$ plays the radius of the sphere.
\item\qquad $\bullet$
The finite matrix representation of $su(2)$. When
$\kappa^2(N^2-1)=4R^2$ then $M_N(\Cmpx)$ forms a quotient algebra, and
$R^2$ is the Casimir operator.
\item\qquad $\bullet$
A noncommutative algebra of polynomials which is 
an infinite dimensional representation of $su(2)$, for other values of
$\kappa$.
\end{trivlist}

In section \ref{ch_Pexp} we introduce the algebra and give a bilinear
form on it. In section \ref{ch_Pnm} we give a basis of
$\Pexpr(\kappa,R)$ which is orthogonal with respect to this bilinear
form. Some of the basis elements were given previously for the matrix
case in \cite{Bayen1}.

In section \ref{ch_Jexp} we give an alternative representation
of the elements of $\Pexpr(\kappa,R)$, and show how the basis elements
can be written in terms of Hahn polynomials.  We also show that
$\Pexpr(\kappa,R)$ may be viewed as an infinite dimensional,
reducible, non-decomposable representation of $su(2)$, and give some
of its ideals.

In section \ref{ch_S2} we look at the commutative case, $\kappa=0$,
and show that the basis elements become the standard spherical
harmonics.  In section \ref{ch_Moy} we calculate the Moyal bracket
which is the limit of the commutator as $\kappa\to0$. We also look at
what the limits of the standard operators on ${\Pexpr(\kappa,R)}$
are.

Finally in the appendix A we draw attention to some facts about
the universal enveloping algebra of $su(2)$ which are needed for some
of the proofs.


\subsection*{Notation and order of proofs}

The summation convention is not used in any part of the article.
The results in the appendix A are used throughout
the article. This section may be read first since no proofs in this
section require material from the rest of the article.


\section{ Definition of the algebra $\Pexpr(\kappa,R)$}
\label{ch_Pexp}

Given the constants $\kappa,R\in\Real$, with $R>0$,
we define the algebra $\Pexpr(\kappa,R)$ to be
\begin{eqnarray}
\Pexpr(\kappa,R) &=&
\{ \mbox{Free noncommuting algebra of polynomials in $x,y,z$ } \}
\Big/\sim
\label{Pex_def_P}
\end{eqnarray}
where $\sim$ are the relations: 
\begin{eqnarray}
[x,y] \sim i\kappa z,\
[y,z] \sim i\kappa x,\
[z,x] \sim i\kappa y,\
x^2+y^2+z^2 \sim R^2
\label{Pex_com_rel}
\end{eqnarray}
We note that this is a well defined algebra since it is equivalent to
the quotient
\begin{eqnarray}
\Sexpr(\kappa) / J(R)
\label{Pex_altdef_P}
\end{eqnarray}
where $\Sexpr=\Sexpr(\kappa)$ is the universal enveloping algebra of
the Lie algebra $su(2)$ and $J(R)$ is the ideal
\begin{eqnarray}
J(R)=\{(x^2+y^2+z^2-R^2)f\ |\ f\in\Sexpr\} 
\label{Pex_def_JR}
\end{eqnarray}
This ideal is two-sided since $(x^2+y^2+z^2-R^2)$ commutes with all
elements in $\Sexpr$.  As shown in section \ref{ch_redu} we may view
$\Pexpr(\kappa,R)$ as a the vector space for a representation of
$su(2)$. This representation is reducible but not decomposable. (The
same of which is true for $\Sexpr(\kappa)$, the universal enveloping
algebra). Any attempt to give $\Pexpr(\kappa,R)$ a Hilbert space
structure would make the this representation of $su(2)$ 
non-continuous.

Usually $\kappa$ and $R$ are implicit and we
simply write $\Pexpr$. We chose the representatives of each equivalent
class $f\in\Pexpr$ to be the totally symmetric formally trace-free
polynomial in $\{x,y,z\}$.  This means that we can write $f\in\Pexpr$
as
\begin{eqnarray}
f &=& \sum_{n=0}^{\degreerm(f)}  \sum_{a_1\ldots a_n=1}^3 
f_{a_1\ldots a_n} x^{a_1}x^{a_2}\ldots x^{a_n}
\label{Pex_f_xxx}
\end{eqnarray}
where $\{x^1,x^2,x^3\}=\{x,y,z\}$. Each $f_{a_1\ldots a_n}$ is
completely symmetric in its indices and satisfies:
\begin{eqnarray}
\sum_{b=1}^3 f_{bba_3a_4\ldots a_n} &=& 0
\label{Pex_form_tr0}
\end{eqnarray}
The condition (\ref{Pex_form_tr0}) will be called formally trace-free
to distinguish it from the matrix trace. In appendix A we
give some more information about the elements of $\Sexpr$ which can be
written as totally symmetric polynomials and we define the formal
trace of these elements.

There is a natural linear bijection 
\begin{eqnarray}
\Psi_{\kappa_1,\kappa_2} &:&
\Pexpr(\kappa_1,R)\mapsto\Pexpr(\kappa_2,R) 
\label{Pex_def_Psi_kk}
\end{eqnarray}
given as follows: Let $f\in\Pexpr(\kappa_1,R)$ and
$g\in\Pexpr(\kappa_1,R)$ are both written in the totally symmetric
formally trace-free form (\ref{Pex_f_xxx}). Then
$\Psi_{\kappa_1,\kappa_2}(f)=g$ if and only if $f_{a_1\ldots
a_n}=g_{a_1\ldots a_n}$ for all indices $\{a_1,a_2,\ldots\}$.  This
mapping is principally used when one of the $\kappa$'s is zero since it
then relates the commutative algebra of functions on the sphere with
the noncommutative algebra. It is clear from the definition of
$\Psi_{\kappa_1,\kappa_2}$ that it is not a homomorphism. (i.e. it
does not preserve the product on $\Pexpr(\kappa,R)$ ).

Let $\Pexpr^n\subset \Pexpr$ be the set of all homogeneous 
polynomials of order $n$, i.e.
\begin{eqnarray}
\Pexpr^n &=& \left\{\sum_{a_1\ldots a_n=1}^3 f_{a_1a_2\cdots a_n}
x^{a_1}x^{a_2}\cdots x^{a_n} \ \bigg|\hbox{ $f_{a_1a_2\cdots a_n}$ is
totally symmetric and formally trace-free}\! \right\}\quad
\label{Pex_def_Pn}
\end{eqnarray}
Then $\dim(\Pexpr^n)=2n+1$. So as a set
\begin{eqnarray}
\Pexpr &=& \bigoplus^\infty_{n=0} \Pexpr^n  
\qquad\hbox{finite sums only}
\end{eqnarray}
We define the projection
\begin{eqnarray}
\pi_n &:& \Pexpr\mapsto \bigoplus_{r=1}^n \Pexpr^r
\end{eqnarray}
We define the operation of taking the Hermitian conjugate by
\begin{eqnarray}
\dagger:\Pexpr\mapsto\Pexpr,\
(ab)^\dagger = b^\dagger a^\dagger,\ 
x^\dagger=x,\ 
y^\dagger=y,\ 
z^\dagger=z,\ 
\lambda^\dagger=\cnj\lambda \qquad\hbox{for $\lambda\in\Cmpx$}
\end{eqnarray}
There is a sesquilinear form on $\Pexpr$ given by
\begin{eqnarray}
\langle\bullet,\bullet\rangle &:& 
\Pexpr\times\Pexpr \mapsto \Cmpx
\cr
\langle f,g \rangle &=&
\pi_0(f^\dagger g)
\label{Pex_def_IP}
\end{eqnarray}
In the section \ref{ch_Pnm} we give a basis of $\Pexpr^n$ and $\Pexpr$
which are orthogonal with respect to this bilinear form.  We also show
that this form is Hermitian $\langle f,g \rangle = \cnj{\langle g,f
\rangle}$. However this bilinear form is not positive definite and
$\langle f,f \rangle$ may be positive, negative or zero. It could be
called a {\it degenerate pseudo inner product}.

As stated in the introduction, there are a number of values of
$\kappa$ and $R$ for which $\Pexpr(\kappa,R)$ is a special algebra. If
$\kappa=0$ then $\Pexpr(0,R) \cong \Czz(S^2)$, the set of finite sums
of spherical harmonics. In this case $\langle\bullet,\bullet\rangle$
does become positive definite and equal to the standard
inner product of functions on $S^2$. We can then close $\Czz(S^2)$ to
give $L^2(S^2)$. This will be analised in section \ref{ch_S2}.

If $\kappa^2(N^2-1)=4R^2$ with $N\in\Intg$, $N\ge1$ then
$\pi_{N-1}\Pexpr(\kappa,R) \cong M_N$ the set of
$N\times N$ matrices.  This isomorphism is given explicitly in section
\ref{ch_Rep}.  In this case the bilinear form restricted to
$\Pexpr^{N-1}$ is an inner product, while the bilinear form on any
other element vanishes.

It is also possible to take the limit $R\to\infty$, with $\kappa$
constant.  This case is dealt with in~\cite{Madore2}.


\section{ Orthogonal Basis of $\Pexpr(\kappa,R)$ }
\label{ch_Pnm}

In this chapter we give an explicit orthonormal basis for $\Pexpr^n$
and $\Pexpr$.  As when dealing with representations on $su(2)$, it
turns out to be very convenient to work with the ladder operators
\begin{eqnarray}
J_+ = x+iy, \qquad J_- = x-iy
\label{Pmn_def_JpJm}
\end{eqnarray}
which obey the relations
\begin{eqnarray}
[z,J_+] = \kappa J_+ ,\ 
[z,J_-] = -\kappa J_- ,\
[J_+,J_-] = 2\kappa z ,\
z^2 + \tfrac12 J_-J_+ + \tfrac12 J_+J_- &=& R^2 
\label{Pmn_rel_JJ}
\end{eqnarray}


\vskip 1em
Useful operators on $\Pexpr$ are given by
\begin{eqnarray}
e_z f &=& [J_z,f] \qquad \hbox{and simiarly for $e_x,e_y$} 
\label{Pmn_def_ez}
\\
e_\pm f &=& [J_\pm,f] 
\label{Pmn_def_epm}
\\
\Delta  &=& e_z^2 - \kappa e_z + e_+ e_- = e_x^2 + e_y^2 + e_z^2 
\label{Pmn_def_Del}
\end{eqnarray}
We will call $\Delta$ the Laplacian operator.

\begin{lemma}
\label{Pnm_lm_ep_dag}
The Laplacian $\Delta$ commutes with $e_+,e_-,e_z$.  With respect to
the bilinear form {\textup (\ref{Pex_def_IP})}, $e_z$ and $\Delta$ are
self adjoint whilst $e_-^\dagger= e_+$.
\end{lemma}

\begin{proof}{}
\begin{eqnarray}
\langle e_- f,g \rangle &=&
\pi_0( (J_-f)^\dagger g - (fJ_-)^\dagger g) 
=
\pi_0( f^\dagger J_+ g - J_+ f^\dagger g)
\cr &=&
\pi_0( f^\dagger J_+ g - f^\dagger g J_+) -\pi_0(e_+(f^\dagger g))
=
\langle f,e_+ g \rangle \qquad\forall f,g\in\Pexpr
\nonumber
\end{eqnarray}
since from corollary \ref{Sex_corol}: $\pi_0(e_+ f)=0$ for all
$f\in\Pexpr$.
\end{proof}
The effects of these operators are calculated below.
We will show that the spaces $\Pexpr^n$ are invariant under the
operators $e_x,e_y,e_z,e_\pm,\Delta$,
and that the $\Pexpr^n$ are the orthogonal eigenspaces of $\Delta$.

These operators vanish when $\kappa=0$ and we obtain the commutative
algebra of functions on the sphere. In section \ref{ch_Moy} we give
the value $\kappa^{-1}e_+$ etc and $\kappa^{-2}\Delta$. The last of
these is the standard Laplacian for functions on the sphere.

\vskip 1em

In the following theorem we give a basis of $\Pexpr$ which is
orthogonal with respect to the bilinear form (\ref{Pex_def_IP}).
This is given by
$P^n_m$ where $n,m\in\Intg,n\ge0,|m|\le n$. The 
$P^n_m$ are defined to be proportional to $e_-^{n-m}(J_+^n)$, and
normalised so that $\langle P^m_n,P^m_n \rangle\in\{1,0,-1\}$.

\begin{theorem}
\label{thm_Pmn}
For $\kappa\ne0$ there is a basis of $\Pexpr(\kappa,R)$ given by 
\begin{eqnarray}
\{P^m_n(\kappa,R)\,|\ n,m\in\Intg,n\ge0,|m|\le n\}
\end{eqnarray}
where
\begin{eqnarray}
P^m_n(\kappa,R) &=&
\alpha_n\kappa^{m-n}
\left(\frac{(n+m)!}{(2n)!\,(n-m)!}\right)^{\scrhalf} 
e_-^{n-m}(J_+{}^n)
\label{Pmn_def_Pmn}
\end{eqnarray}
We shall write $P^m_n=P^m_n(\kappa,R)$ when there is no doubt about
$\kappa$ and $R$. This basis is orthogonal with respect to the
bilinear form.  The ``normalisation'' constant $\alpha_n\in\Real,\
\alpha_n>0$ defined so that $\langle P^m_n,P^m_n
\rangle\in\{1,0,-1\}$.

Each $P^m_n$ can be written as a homogeneous formally trace-free
symmetric polynomial in $(x,y,z)$ of order $n$. Thus the set
$\{P^m_n,m=-n\ldots n\}$ forms an orthogonal basis for $\Pexpr^n$ and
\begin{eqnarray}
\Psi_{\kappa_1,\kappa_2}(P^m_n(\kappa_1,R)) &=& 
P^m_n(\kappa_2,R)
\end{eqnarray}

Each $P^m_n$ is an eigenvector of the operators $e_z$
and $\Delta$.
\begin{eqnarray}
e_z P^m_n &=& \kappa m P^m_n 
\label{Pmn_ez_Pmn} 
\\
\Delta P^m_n &=& \kappa^2 n(n+1) P^m_n 
\label{Pmn_Del_Pmn} 
\end{eqnarray}
so that the $\Pexpr^n$ are the orthogonal eigenspaces of $\Delta$.  The
ladder operators $e_+,e_-$ increase or decrease $m$ so that $\Pexpr^n$
can be viewed as a $2n+1$ dimensional adjoint representation of
$su(2)$.
\begin{eqnarray}
e_+ P^m_n &=& \kappa (n-m)^{\scrhalf} (n+m+1)^{\scrhalf} P^{m+1}_n 
\label{Pmn_ep_Pmn}
\\
e_- P^m_n &=& \kappa (n+m)^{\scrhalf} (n-m+1)^{\scrhalf} P^{m-1}_n 
\label{Pmn_em_Pmn}
\end{eqnarray}

The effect of taking the Hermitian conjugate is given by
\begin{eqnarray}
(P^m_n)^\dagger &=& (-1)^m P^{-m}_n 
\label{Pmn_Pmn_dag}
\end{eqnarray}

\end{theorem}


\begin{proof}{}
To show $P^n_m$ is an eigenvector of $e_z$ (\ref{Pmn_ez_Pmn}) we have
\begin{eqnarray}
e_z P^n_n &=&
\alpha_n e_z J_+{}^n 
=
\alpha_n \sum_{r=0}^{n-1} J_+{}^r e_z(J_+) J_+{}^{n-r-1}
=
\kappa n P^n_n
\nonumber 
\end{eqnarray}
whilst
\begin{eqnarray}
e_z e_- P^m_n &=& e_- e_z P^m_n + \ad_{[z,J_-]} P^m_n 
=
\kappa m \, e_- P^m_n - \kappa  e_- P^m_n
=
\kappa(m-1) e_- P^m_n
\nonumber  
\end{eqnarray}
Thus (\ref{Pmn_ez_Pmn}) follows by induction.


For the ladder operators $e_+,e_-$,
(\ref{Pmn_em_Pmn}) is obvious from the definition of $P^m_n$. To show
(\ref{Pmn_ep_Pmn}) we note
\begin{eqnarray}
e_+(e_-^{n-m} J_+^n) &=&
[e_+ ,e_-^{n-m}] J_+^n 
=
2\kappa \sum_{r=1}^{n-m} e_-^{r-1} e_z e_-^{n-m-r} J_+^n 
\cr &=&
2\kappa^2 \sum_{r=0}^{n-m-1} e_-^{r-1} e_-^{n-m-r}  J_+^n (m+r)
=
2\kappa^2 e_-^{n-m-1} J_+^n \sum_{r=0}^{n-m-1}(m+r)
\cr &=&
\kappa^2 (n-m)(n+m+1) e_-^{n-m-1} J_+^n 
\nonumber
\end{eqnarray}
From the effects of $e_+,e_-,e_z$ (\ref{Pmn_Del_Pmn}) is trivial. 
Orthogonality is simply an application of $e_z$, $\Delta$ and lemma
\ref{Pnm_lm_ep_dag}.   


To show that $\alpha_n$ is independent of $m$, we note
\begin{eqnarray}
\langle P^m_n, P^m_n \rangle &=&
\kappa^{-1}(n+m+1)^{-\scrhalf}(n-m)^{-\scrhalf}
\langle e_- P^{m+1}_n,P^m_n \rangle
\cr
&=&
\kappa^{-1}(n+m+1)^{-\scrhalf}(n-m)^{-\scrhalf}
\langle P^{m+1}_n, e_+ P^m_n \rangle
\cr
&=&
\langle P^{m+1}_n, P^{m+1}_n \rangle
\nonumber
\end{eqnarray}


For the effect of Hermitian conjugation we see that since
$e_z(P^0_n)=0$ then $P^0_n$ may be written as a polynomial in $z$, and
from (\ref{Pmn_rel_JJ}) it is real.  This polynomial given explicitly
in section \ref{ch_Jexp}.  Therefore $(P^0_n)^\dagger=P^0_n$. Since
\begin{eqnarray}
(e_+ f)^\dagger =
[J_+,f]^\dagger =
[f^\dagger,J_-] = -e_-(f^\dagger)
\nonumber
\end{eqnarray}
then
\begin{eqnarray}
(e_+^m f)^\dagger = (-1)^m e_-^m(f^\dagger)
\nonumber
\end{eqnarray}
By repeated use of (\ref{Pmn_ep_Pmn}) and (\ref{Pmn_em_Pmn}) we have
\begin{eqnarray}
P^m_n &=& \kappa^{-m} \left(\frac{(n+m)!}{(n-m)!} \right)^\scrhalf
e_+^m(P^0_n)
\cr
P^{-m}_n &=& \kappa^{-m} \left(\frac{(n+m)!}{(n-m)!} \right)^\scrhalf
e_-^m(P^0_n)
\nonumber
\end{eqnarray}
Thus (\ref{Pmn_Pmn_dag}) follows.


The proof of the statement that $\{P^m_n\,|\, m=-n\ldots n\}$ forms a
basis for $\Pexpr^n$ is given in appendix A.
\end{proof}

We note that the elements $P^0_n$ were discovered by Bayen and Fronsdal
\cite{Bayen1}. However they don't mention using the ladder operators
to get all the elements.


\subsection*{Finite Representation of $\Pexpr$}
\label{ch_Rep}

\begin{theorem}
\label{thm_Rep_vphiN}
For the discrete set of $\kappa$ such that 
\begin{eqnarray}
\kappa^2(N^2-1)=4R^2
\label{Rep_NRk}
\end{eqnarray}
where $N\in\Intg$, $N\ge1$ there exist a surjective homomorphism 
\begin{eqnarray}
\varphi_N &:& \Pexpr(\kappa,R) \mapsto M_N(\Cmpx)
\end{eqnarray}
This is the $N\times N$ representation of $su(2)$. It satisfies
\begin{eqnarray}
\varphi_N(fg) &=& \varphi_N(f)\varphi_N(g) 
\qquad \forall f,g\in\Pexpr
\label{Rep_homo}
\\ 
\varphi_N(f)&=&0 
\qquad \forall f\in\Pexpr^n,\ n\ge N 
\label{Rep_Nf0}
\\
\pi_0(f) &=& \tfrac1N\tr(\varphi_N(f)) 
\qquad\forall f\in\Pexpr
\label{Rep_pi0tr}
\end{eqnarray}
The bilinear form is a genuine inner product if restricted to
$\pi_{N-1}(\Pexpr)$, whilst $\langle f,g \rangle=0$ if $f\in\Pexpr^n$,
$g\in\Pexpr$ and $n\ge N$.
\end{theorem}

\begin{proof}{}
If $\varphi_N$ is an $M_N(\Cmpx)$ representation of $su(2)$ then
$\varphi_N(x),\varphi_N(y),\varphi_N(z)$ satisfy the same commutation
relations, and the Casimir operator is the same as in
$\Pexpr(\kappa,R)$. Thus $\varphi_N$ is a homomorphism
(\ref{Rep_homo}).

As in the example below one can write $\varphi_N(J_+)$ as an upper
triangular matrix (with zeros on the diagonal). Therefore
$\varphi_N(J_+^N)=0$. Since $\varphi_N$ is a homomorphism
$\varphi_N(e_-(f))=0$ if $\varphi_N(f)=0$, so $\varphi_N(P^m_n)=0$ if
$n\ge N$. Hence (\ref{Rep_Nf0}).

Since can write $\varphi_N(J_+)$ as an upper
triangular matrix $\tr(\varphi_N(J_+^n))=0$ for all $n>0$, and
$\tr(\varphi_N(e_+(f)))=\tr[\varphi_N(J_-),\varphi_N(f)]=0$ so 
$\tr(\varphi_N(P^n_m))=0$ for all $n>0$. This gives (\ref{Rep_pi0tr}).
\end{proof}


We now give an explicit representation $\varphi_N:\Pexpr\mapsto
M_N(\Cmpx)$.  This representation is very useful for calculating
formulae.  Let $\vvec(N,r)$ with $r=0\ldots N-1$ be the orthogonal
basis column vectors which are eigenvectors of $\varphi_{N}(z)$. Since
the dimension of the representation is explicit we drop $\varphi_N$.
Then by rewriting the standard ladder operators we have
\begin{eqnarray}
J_+ \vvec(N,r) &=& \kappa (N-r-1)^{\scrhalf} (r+1)^{\scrhalf} 
\vvec(N,r+1) 
\label{Rep_Jp_vec}
\\
J_- \vvec(N,r) &=& \kappa r^{\scrhalf} (N-r)^{\scrhalf} \vvec(N,r-1) 
\label{Rep_Jm_vec}
\\
z \vvec(N,r) &=& \kappa (r-\tfrac{N-1}2) \vvec(N,r) 
\label{Rep_z_vec}
\end{eqnarray}
The dual to $\vvec(N,r)$ is given by $\lvec(r,N)$.
The usefulness of these representation is given by the following
theorem.

\begin{theorem}
\label{thm_Rep}
Let $f$ be a (noncommuting) polynomial in $\{J_+,J_-,z,\kappa,R\}$.
Using the identities (\ref{Pmn_rel_JJ}) we can put $f\in\Pexpr$.
Thus we can write 
\begin{eqnarray}
f=\sum_{n=0}^{\degreerm(f)} \sum_{m=-n}^n 
f_{nm}(\kappa,R) P^m_n
\end{eqnarray}
where the product $\alpha_n(\kappa,R)f_{nm}(\kappa,R)$ is a polynomial
in $\kappa$ and $R$. The $f_{nm}$ can be calculated simply from its
values when $\kappa^2(N^2-1)=4R^2$. That is if $g_{nm}(\kappa,R)$ is
another polynomial and
\begin{eqnarray}
f_{nm}(2R(N^2-1)^{-\scrhalf},R) &=& g_{nm}(2R(N^2-1)^{-\scrhalf},R)
\qquad \forall N\in\Intg, N\ge 2
\end{eqnarray}
Then $f_{nm}(\kappa,R)=g_{nm}(\kappa,R)$ for all $\kappa,R$.
This result is independent of the value of $\alpha_n(\kappa,R)$.
\end{theorem}

\begin{proof}{}
The basis term $P^m_n/\alpha_n$ is a polynomial in
$\{J_+,J_-,z,\kappa,R\}$. Its value is independent of the actual
definition of $\alpha_n$. The manipulation of $f$ into the above form
makes sure that $\alpha_n f_{nm}$ is a polynomial.  The second part
follows since all polynomials are determined by there value on a
finite number of distinct points.
\end{proof}

We note that given a polynomial $f(J_+,J_-,z)$ of order $r$. We can
use the above theorem to calculate $f_{nm}$. This process is of the
order of $r^3$.  However if we directly use the equations
(\ref{Pmn_rel_JJ}) then the process takes an exponential amount of
time. This can be used to get computers to calculate explicit
expressions in the $P^m_n$.

We can now extend some of the basic facts about the matrix trace for
all $\kappa,R\in\Real$ and $R>0$
\begin{corol}
\label{Rep_sym_pi0}
\begin{eqnarray}
\pi_0(fg) &=& \pi_0(gf) \textup{ and } 
\langle f,g \rangle = \cnj{\langle g,f \rangle}
\qquad\forall f,g\in\Pexpr,\quad \kappa,R\in\Real, R>0
\end{eqnarray}
\end{corol}

\subsection*{The ``normalisation'' constant $\alpha_n$}

\begin{theorem}
\label{thm_Norm_alpha_n}
Since
\begin{eqnarray}
\pi_0(J^n_-J^n_+) &=&
\frac{(n!)^2}{(2n+1)!} 
\prod_{r=1}^n(4R^2+\kappa^2(1-r^2))
\label{Norm_pi0_JnJn}
\end{eqnarray}
we define
\begin{eqnarray}
\sigma_n(\kappa,R) &=& \sign(\pi_0(J^n_-J^n_+)) 
\label{Norm_def_sig_n}
\\
\alpha_n(\kappa,R) &=&
\left\{
\begin{array}{ll}
|\pi_0(J^n_-J^n_+)|^{-\scrhalf} 
&\qquad
\sigma_n(\kappa,R) \ne 0
\cr
1 &\qquad
\sigma_n(\kappa,R) = 0
\end{array}
\right.
\label{Norm_def_alpha_n}
\end{eqnarray}
If we let $N_0$ be the smallest integer greater than or equal to
$(4R^2\kappa^{-2}+1)^\scrhalf$. i.e.
\begin{eqnarray}
N_0 &=& \lceil (4R^2\kappa^{-2}+1)^\scrhalf \rceil
\end{eqnarray}
Then the value of $\sigma_n(\kappa,R)$ is given by
\begin{eqnarray}
\sigma_n(\kappa,R) &=&
\left\{
\begin{array}{ll}
1 & n\le N_0-1 
\cr
(-1)^{n-N_0+1} \qquad 
& n\ge N_0,\textup{ and } (4R^2\kappa^{-2}+1)^\scrhalf
\not\in\Intg 
\cr
0 &
n\ge N_0, \textup{ and } (4R^2\kappa^{-2}+1)^\scrhalf=N_0\in\Intg 
\end{array}
\right.
\end{eqnarray}
The ``normalisation'' of $P^m_n$ is now
\begin{eqnarray}
\langle P^m_n, P^m_n \rangle &=& \sigma_n(\kappa,R)
\end{eqnarray}
\end{theorem}

\begin{proof}{}
By repeated application of (\ref{Rep_Jm_vec}) we have
\begin{eqnarray}
\lvec(r,N)J_-^{m}J_+^{m}\vvec(N,r)
&=&
\kappa^{2m}
\frac{(N-r-1)!}{(N-r-m-1)!}
\frac{(r+m)!}{r!}
\nonumber
\end{eqnarray}
thus 
\begin{eqnarray}
\pi_0(\varphi_N(J_-^mJ_+^m))
&=&
\kappa^{2m}\frac1N
\sum_{r=0}^{N-1}
\frac{(N-r-1)!}{(N-r-m-1)!}
\frac{(r+m)!}{r!}
\cr
&=&
\kappa^{2m}
\frac{m!(N-1)!}{N(N-m-1)!}
\sum_{r=0}^{N-m-1}
\frac{(1+m-N)_r (m+1)_r}{(1-N)_r r!}
\cr
&=&
\kappa^{2m}
\frac{m!(N-1)!}{N(N-m-1)!}
F(1+m-N,m+1;1-N;1)
\nonumber
\end{eqnarray}
where $F(a,b;c;z)$ is the hypergeometric function. Since this has only
a finite number of terms we may write
\begin{eqnarray}
F(1+m-N,m+1;1-N;1) &=&
\lim_{\varepsilon\to 0} 
F(1+m-N,m+1;1-N+\varepsilon;1)
\nonumber
\end{eqnarray}
which is also a polynomial. We may then use the standard result
\begin{eqnarray}
F(a,b;c;1)  &=& \frac{\Gamma(c)\Gamma(c-a-b)}{\Gamma(c-a)\Gamma(c-b)}
\nonumber
\end{eqnarray}
This formula is normally only valid when $\Re(c-a-b)>0$, however, here
it is valid here since the number of terms is finite.  Thus
\begin{eqnarray}
F(1+m-N,m+1;1-N;1) &=&
\lim_{\varepsilon\to 0} 
\frac{\Gamma(1-N+\varepsilon) \Gamma(\varepsilon-2m-1)}%
{\Gamma(-m+\varepsilon)\Gamma(-N-m+\varepsilon)}
=
\frac{m!(N+m)!}{(2m+1)!(N-1)!}
\nonumber
\end{eqnarray}
Which implies
\begin{eqnarray}
\pi_0(J_-^mJ_+^m) &=&
\kappa^{2n} 
\frac{(N+n)!}{N(N-n-1)!}
\frac{(n!)^2}{(2n+1)!} 
\nonumber
\end{eqnarray}
substituting $N=(4R^2\kappa^{-2}+1)^\scrhalf$ and using theorem
\ref{thm_Rep} gives (\ref{Norm_pi0_JnJn}). The rest of the theorem is
derived from this equation.
\end{proof}



\section{``$\Jexpr$-expressions'', an Alternative Representations for 
 $\Pexpr(\kappa,R)$}
\label{ch_Jexp}

Given $f\in\Pexpr$ then we can use the commutation relations
(\ref{Pmn_rel_JJ}) to push the $J_+$ and $J_-$ to the left of each
term. If a $J_+$ and $J_-$ appear in one term we can use the Casimir
identity to remove both. Thus the resulting terms must either have
only $J_+$'s or only $J_-$'s or neither.  If we collect all the terms
with the same number of $J_+$ or $J_-$ as there factors then $f$ we
may written as a sum of terms of the form
\begin{eqnarray}
\{(J_+)^m p(z)\,,\, p(z)\,,\, (J_-)^{-m}p(z)\}
\nonumber
\end{eqnarray}
where $p(z)$ is a polynomial in $z$. By looking at the action of $e_z$
on each of these possibilities we see that if $f$ is an eigenvector of
$e_z$ then we my write
\begin{eqnarray}
f = \left\{
\begin{array}{ll}
(J_+)^m p^m_f(z,\kappa,R)    & 
\hbox{ if $e_zf = \kappa m f$ and $m>0$ } \cr
p^0_f(z,\kappa,R)            & \hbox{ if $e_zf = 0$ } \cr
(J_-)^{-m} p^m_f(z,\kappa,R) & 
\hbox{ if $e_zf = \kappa m f$ and $m<0$ } \cr
\end{array}
\right.
\label{J_f_Jp}
\end{eqnarray}
This will be known as $\Jexpr$ notation. It turns our to be more
convenient to consider only the case when $\Pexpr$ has a finite
representation, and we therefore write $p^m_f(z,N)$. This is a
function of $z$, $N$ and $\kappa$ with $\kappa$ being implicit.
However all the expressions can be generalised by substituting
$N=(4R^2\kappa^{-2}+1)^\scrhalf$, and using theorem \ref{thm_Rep}.
This theorem is valid since the definition of $\kappa^n
P^m_n/\alpha_n$ automatically makes it a polynomial in
$(J_\pm,z,R,\kappa)$

\begin{theorem}
\label{thm_J}
Since $P^m_n$ is a eigenvector of $e_z$ we have
\begin{eqnarray}
P^m_n = \left\{
\begin{array}{ll}
(J_+)^m p^m_n(z,N)    & \hbox{\rm\ for }  m>0 \cr
p^0_n(z,N)            & \hbox{\rm\ for }  m=0 \cr
(J_-)^{-m} p^m_n(z,N) & \hbox{\rm\ for }  m<0 \cr
\end{array}
\right.
\label{J_Pmn_Jp}
\end{eqnarray}
where for $m\ge 0$
\begin{eqnarray}
p^m_n(z,N) &=&
\alpha_n (-\kappa)^{n-m} \ppmatrix{2n \cr n-m}^{-\scrhalf}
h^{(m,m)}_{n-m}(\nfrac{z}{\kappa}+\tfrac{N-1}2,N-m)
\label{J_thm_Pmn_h}
\\
p^{-m}_n(z,N) &=&
\alpha_n (-1)^m(\kappa)^{n-m} \ppmatrix{2n \cr n-m}^{-\scrhalf}
h^{(m,m)}_{n-m}(\nfrac{z}{\kappa}-m+\tfrac{N-1}2,N-m)
\label{J_thm_Pmmn_h}
\end{eqnarray}
where $h^{(\alpha,\beta)}_n(x,N)$ follows the notation of
\cite[chapter 2]{Nik}. These are the Hahn Polynomials.  This has an
explicit formulation in terms of generalised hypergeometric functions:
\begin{eqnarray}
p^m_n(z,N) &=& 
\alpha_n (-\kappa)^{n-m} \ppmatrix{2n \cr n-m}^{-\scrhalf} 
\frac{(m+1)_{n-m} (m+1-N)_{n-m}}{(n-m)!} \times \cr
&& \qquad
{}_3F_2 (m-n,-\nfrac{z}{\kappa}-\tfrac{N-1}2,n+m+1;m+1,m+1-N;1)  
\label{J_thm_Pmn_F}
\end{eqnarray}
where $m\ge 0$.
\end{theorem}


Before proving this, let's start with a little lemma.
\begin{lemma}
\label{J_lm_Jp}
For a polynomial $p(z)$ and $m\in\Intg^+$ a positive integer
\begin{eqnarray}
p(z) J_+^m  &=& J_+^m p(z+m\kappa) 
\\
p(z) J_-^m  &=& J_-^m p(z-m\kappa) 
\\
J_+^m J_-^m &=& 
\prod_{s=0}^{m-1}
\Big( R^2 - (z-s\kappa)(z-(s+1)\kappa) \Big) 
\\
J_-^m J_+^m &=& 
\prod_{s=0}^{m-1}
\Big( R^2 - (z+s\kappa)(z+(s+1)\kappa) \Big)
\end{eqnarray}
\end{lemma}

\begin{proof}{}
From (\ref{Pmn_rel_JJ}) we have $zJ_+=J_+(z+\kappa)$. Thus
$zJ_+^m=J_+^m(z+m\kappa)$ and $z^pJ_+^m=J_+^m(z+m\kappa)^p$, hence
result.

For the other equations we note 
\begin{eqnarray}
J_+^m J_-^m = J_+^{m-1}(R^2-z(z-\kappa)) J_-^{m-1}
= J_+^{m-1}J_-^{m-1}(R^2-(z-(m-1)\kappa)(z-m\kappa))
\nonumber
\end{eqnarray}
\end{proof}


\begin{proof}{of theorem \ref{thm_J}}
Given two basis harmonics $P^{m_1}_{n_1},P^{m_2}_{n_2}\in\Pexpr$ then
we have $\langle P^{m_1}_{n_1},P^{m_2}_{n_2} \rangle = 0$ if $m_1\ne
m_2$ or $n_1\ne n_2$.  Writing these in $\Jexpr$ form we have for
$m_1\ne m_2$ or $n_1\ne n_2$ and $m_1,m_2\ge 0$
\begin{eqnarray}
\pi_0
\big(p^{m_1}_{n_1}(z,N)\, J_-^{m_1}\, J_+^{m_2} \, p^{m_2}_{n_2}(z,N)
\big)
&=& 0
\nonumber
\end{eqnarray}
thus
\begin{eqnarray}
\sum_{r=0}^{N-1}
\lvec(r,N) p^{m_1}_{n_1}(z,N) \, J_-^{m_1}\, 
J_+^{m_2} \, p^{m_2}_{n_2}(z,N) \vvec(N,r) &=& 0
\nonumber
\end{eqnarray}
It is clear this is satisfied for $m_1\ne m_2$. It is also obvious
that the summand vanished if $r+m>N-1$
Taking $m_1=m_2=m\ge 0$ and
$n_1\ne n_2$ we have
\begin{eqnarray}
\sum_{r=0}^{N-m-1}
p^{m}_{n_1}(\kappa(r-\tfrac{N-1}2),N) \, 
p^{m}_{n_2}(\kappa(r-\tfrac{N-1}2),N) \,
\lvec(r,N)J_-^{m}J_+^{m}\vvec(N,r)
\nonumber
\end{eqnarray}
Now 
\begin{eqnarray}
\lvec(r,N)J_-^{m}J_+^{m}\vvec(N,r)
&=&
\kappa^{2m}\left(
\frac{(N-r-1)!}{(N-r-m-1)!}
\frac{(r+m)!}{r!} \right)
\label{J_Jpm_Vec}
\nonumber
\end{eqnarray}
is precisely the weight function for the Hahn
Polynomials $h^{(\alpha,\beta)}_{n'}(r,N')$ where $\alpha=m$,
$\beta=m$, $n'=n-m$ and $N'=N-m$.  Now to get the correspondence
between the functions we look at coefficient of the highest order in
$r$. It is easy to show
\begin{eqnarray}
e_-(J_+^m z^p) &=& -\kappa (p+2m) J_+^{m-1} ( z^{p+1} + O(z,p) ) 
\nonumber
\end{eqnarray}
where $O(z,p)$ is a polynomial in $z$ of order $p$ or less.
So
\begin{eqnarray}
e_-^p(J_+^m) = (-\kappa)^p \frac{(2n)!}{(2n-p)!} J_+^{n-p} (z^p +
O(z,p-1))
\nonumber
\end{eqnarray}
so
\begin{eqnarray}
P^m_n &=& \alpha_n(-1)^{n-m} \ppmatrix{2n \cr n-m}^{\scrhalf} J_+^m (z^{n-m}
+ O(z,n-m-1))
\nonumber
\end{eqnarray}
so $p^m_n(z,N)$ is a Polynomial in $z$ of order $n-m$. This is the
same order $h^{(m,m)}_{n-m}(r,N)$. From \cite[page 42]{Nik} we have
\begin{eqnarray}
h^{(m,m)}_{n-m}(\nfrac{z}{\kappa} + \tfrac{N-1}2,N) &=&
\ppmatrix{2n \cr n-m} (z^{n-m} + O(z,n-m-1))
\nonumber
\end{eqnarray}
hence (\ref{J_thm_Pmn_h}). Expression
(\ref{J_thm_Pmn_F}) follows from the literature.

Finally (\ref{J_thm_Pmmn_h}) is simply an application of lemma
\ref{J_lm_Jp}.
\end{proof}



\subsection*{The reducibility  and ideals of $\Pexpr(\kappa,R)$}
\label{ch_redu}

\begin{lemma}
\label{redu_lm_omega}
For any $f\in \Pexpr$ we define 
\begin{eqnarray}
&&\omega_n : \Pexpr \mapsto \Pexpr \cr 
&&\omega_n(f) = \sum_{m=-n}^n (P^{m}_{n})^\dagger f P^{m}_{n}
\label{redu_def_omega_n}
\end{eqnarray}
then $\omega_n$ is a self-adjoint operator on $\Pexpr$ with respect to
the bilinear form. It commutes with the operators
$e_+,e_-,e_z,\Delta$.  It is diagonal with respect to the basis
elements $P^m_n$, and the eigenvalues depend only on $n$ so we can
write
\begin{eqnarray}
\omega_n(f) = \omega_{na} f \qquad \forall f\in\Pexpr^a
\label{redu_omega_na}
\end{eqnarray}
where $\omega_{na}$ is a real polynomial of $(n,\kappa)$. 
\begin{eqnarray}
\sum_{m=-n}^n (P^{m}_{n})^\dagger P^{m}_{n}
&=& 
\omega_n(1) 
=
\omega_{n0} 
= 
\sigma_n(\kappa,R)(2n+1)
\label{redu_omega_n1}
\end{eqnarray}
\end{lemma}

\begin{proof}{}
Self-adjointness follows from the definition of $\omega_n$ and
corollary \ref{Rep_sym_pi0}. Whilst the fact that it commutes with
$e_\pm,e_z,\Delta$ follows from direct substitution.
Since $\omega_p(P^m_n)$ is a polynomial then by the operations of
$e_z$ and $\Delta$ it is clear that is must be proportional to
$P^m_n$. Furthermore since $\omega_p$ commutes with $e_-$, then the
eigenvector for the space spanned by $P^m_n$ must be independent of $m$.
Hence (\ref{redu_omega_na}). Since $\omega_p(1)\in\Pexpr^0$ then
(\ref{redu_omega_n1}) follows by considering $\pi_0(\omega_p(1))$.
\end{proof}

\begin{lemma}
\label{redu_lm_ideal}
For all $\kappa,R\in\Real$, $\Pexpr(\kappa,R)$ has at least one proper
left ideal given by $I= \{fz\ |\ f\in\Pexpr \}$. Also
$\Pexpr(\kappa,R)$ has a proper two sided idea if and only if
$(4R^2\kappa^{-2}+1)^{\scrhalf}\in\Intg$
\end{lemma}

\begin{proof}{}
To see that $I$ is a proper left ideal we note that $1\not\in I$. For
assume there exists $f\in\Pexpr$ such that $fz=1$, then writing $f$ in
$\Jexpr$-notation we have
\begin{eqnarray}
fz &=& 
\sum_{a=0}^{\degreerm(f)} J_+^a p_a(z)z + 
\sum_{a=1}^{\degreerm(f)} J_-^a p_{-a}(z)z 
=1 
\nonumber
\end{eqnarray}
where $p_a(z)$ is a polynomial in $z$. By looking at the operation of
$e_z$ implies $p_a(z)=0$ for $a\ne 0$, whilst $p_0(z)z=1$ which is
impossible. Thus $I$ is a left ideal of $\Pexpr$. 

We note however that if $\Pexpr$ contained infinite (unbounded) sums
then there is a solution $f=\sum_{n=0}^\infty f_n P^0_n$ such that
$fz=1$. This expression cannot be written in $\Jexpr$-notation.

If $4R^2=\kappa^2(N^2-1)$, for some $N\in\Intg$, then from theorem
\ref{thm_Rep_vphiN} the subspace $\oplus_{r=N}^\infty
\Pexpr^r\subset\Pexpr(\kappa,R)$ is a two sided ideal. This is because
it is the kernel of $\varphi_N$.  Otherwise let $I\subset\Pexpr$ be a
two sided ideal and $f\in I$. Then by the operation of $e_z$ and
$\Delta$ we can show there is a basis element $P^m_n\in I$. By
application of $e_\pm$ we have $\Pexpr^n\subset I$ for some $n$.  From
(\ref{redu_omega_n1}) we have $\sigma_n(\kappa,R)(2n+1)\in I$. Since
$(4R^2\kappa^{-2}+1)^{\scrhalf}\not\in\Intg$ we have from theorem
\ref{thm_Norm_alpha_n}, $\sigma_n(\kappa,R)\ne0$ so $1\in I$.

\end{proof}

\begin{theorem}
The map given by
\begin{eqnarray}
&&\rho:su(2)\mapsto 
\{ f:\Pexpr\mapsto\Pexpr \ | 
\mbox{\textup{ $f$ is linear}} \} \cr
&&\rho(a)(f) = af\qquad a\in su(2),f\in\Pexpr
\end{eqnarray}
may be viewed as an infinite dimensional representation of
$su(2)$. This representation is always reducible but not decomposable.
It is ``Hermitian'' in that it respects the Hermitian conjugate with
defined by the bilinear form
\begin{eqnarray}
\langle \rho(a)f,g\rangle &=& 
\langle f,\rho(a^\dagger) g\rangle
\qquad\forall a\in su(2),\ f,g\in\Pexpr
\end{eqnarray}

\end{theorem}

\begin{proof}{}
The subspace $I\in\Pexpr$ given in lemma \ref{redu_lm_ideal} is
invariant under the action of $\rho$. However $\Pexpr$ is not
decomposable because from the action of $\rho$ on the element
$1\in\Pexpr$ one can generate $\Pexpr$. The Hermitian conjugate is by
direct substitution. 
\end{proof}

We note that the universal enveloping, $\Sexpr(\kappa)$ is also
reducible but not decomposable. 
Any attempt to give $\Pexpr(\kappa,R)$ a Hilbert space structure 
would mean the action of $su(2)$ were non-continuous operators.



\section{The Commutative Case $\kappa=0$}
\label{ch_S2}

As mentioned in the introduction the algebras $\Pexpr(\kappa,R)$ is a
commutative algebra when $\kappa=0$ and is isomorphic to the algebra of
functions on the sphere. In this chapter we will show that
$\Pexpr(0,R)=\Czz(S^2)$ the set of finite sums of spherical harmonics,

To make this isomorphism explicit we write
\begin{eqnarray}
x|_{\kappa=0} =
\fs{x} &=& R\sin\phi\sin\theta \cr
y|_{\kappa=0} =
\fs{y} &=& R\cos\phi\sin\theta \cr
z|_{\kappa=0} =
\fs{z} &=& R\cos\theta 
\label{lim_xyz}
\end{eqnarray}
To distinguish elements of $\Pexpr(\kappa,R)$ with $\kappa\ne0$ from
the elements of $\Czz(S^2)$ the latter are written in bold when there
may be doubt. From (\ref{Pmn_def_JpJm}) we have
\begin{eqnarray}
\fs{J}_+ &=& i e^{-i\phi} R\sin\theta  \cr
\fs{J}_- &=& -i e^{i\phi}  R\sin\theta
\label{lim_Jpm}
\end{eqnarray}
From (\ref{Norm_def_alpha_n}) we have
\begin{eqnarray}
\alpha_n|_{\kappa=0} &=& \frac{((2n+1)!)^{\scrhalf}}{n!} 
(2R)^{-n}
\end{eqnarray}
From (\ref{Rep_NRk}) we may think of the case $\kappa=0$ as the limit
as $N\to\infty$. Near this limit (i.e. for large $N$)
\begin{eqnarray}
\kappa\sim 2R/N
\end{eqnarray}
The definition for $P^m_n$ (\ref{Pmn_def_Pmn})
is not valid in the case $\kappa=0$. We therefore define them as the
limit
\begin{eqnarray}
P^m_n(0,R) &=& \lim_{\kappa\to0} \Psi_{\kappa,0}(P^m_n(\kappa,R))
\label{lim_def_Pmn}
\end{eqnarray}

\begin{theorem}
\label{lim_thm}
In the case $\kappa=0$ the ``fuzzy'' spherical harmonics become the
standard spherical harmonics
\begin{eqnarray}
P^m_n|_{\kappa=0}
&=&
(-1)^n 
\left(\frac{(n+m)!(2n+1)}{(n-m)!}\right)^{\scrhalf}
e^{-im\phi} P^{-m}_n(\cos\theta) 
=
(-1)^n Y^{-m}_n(\theta,\phi) 
\label{lim_Pmn}
\end{eqnarray}
where $P^{m}_n(\nfrac{\fs{z}}{R})$ are the Associated Legendre
functions, and $Y^{-m}_n(\theta,\phi)$ are the orthonormal harmonics
on the sphere.  So $\Pexpr(0,R)=\Czz(S^2)$, the set of finite sums of
spherical harmonics.  The bilinear form on $\Pexpr$ becomes the
standard inner product on $\Czz(S^2)$:
\begin{eqnarray}
\langle f,g\rangle \to \langle \fs{f},\fs{g}\rangle_{S^2}
=\frac{1}{4\pi R^2}\int_{S^2} \cnj{\fs{f}}\fs{g}\sin\theta d\phi
d\theta
\label{lim_def_IP}
\end{eqnarray}

\end{theorem}


\begin{proof}{}
From (\ref{J_thm_Pmn_h}),(\ref{lim_def_Pmn}) and theorem
\ref{thm_Rep}, we have for $m\ge 0$
\begin{eqnarray}
P^m_n|_{\kappa=0} &=& 
\lim_{N\to\infty} (P^m_n) 
\cr
&=&
\lim_{N\to\infty} 
\left( J_+^m 
\alpha_n (-\kappa)^{n-m} \ppmatrix{2n \cr n-m}^{-\scrhalf}
h^{(m,m)}_{n-m}(\nfrac{z}{\kappa}+\tfrac{N-1}2,N-m)
\right)
\cr
&=&
(iR\sin\theta)^m e^{-im\phi}
\frac{((2n+1)!)^{\scrhalf}}{n!} (2R)^{-n} (-2R)^{n-m}
\ppmatrix{2n \cr n-m}^{-\scrhalf} \times
\cr
&&\qquad\qquad
\lim_{N\to\infty} 
\left( N^{m-n}
h^{(m,m)}_{n-m}(\nfrac{z}{\kappa}+\tfrac{N-1}2,N-m)
\right)
\nonumber
\end{eqnarray}
From \cite[page 46]{Nik} this is given by
\begin{eqnarray}
\lim_{N\to\infty} (P^m_n) &=&
(iR\sin\theta)^m e^{-im\phi}
\frac{((2n+1)!)^{\scrhalf}}{n!} (2R)^{-n} (-2R)^{n-m}
\ppmatrix{2n \cr n-m}^{-\scrhalf}
P^{(m,m)}_{n-m}(\cos\theta)
\nonumber
\end{eqnarray}
where $P^{(m,m)}_{n-m}(\nfrac{\fs{z}}{R})$ is the Jacobi Polynomial.
This is related to the Associated Legendre functions by
\begin{eqnarray}
P^{(m,m)}_{n-m}(\cos\theta) &=& \frac{n!}{(n-m)!}
(2i\sin\theta)^{-m} P^{-m}_n(\cos\theta)
\nonumber
\end{eqnarray}
Hence (\ref{lim_Pmn}). For $m>0$ we note that taking the limit of
(\ref{Pmn_Pmn_dag}) 
\begin{eqnarray}
P^{-m}_n|_{\kappa=0} &=& 
(-1)^m \cnj{P^{m}_n}|_{\kappa=0} =
(-1)^n Y^m_n(\theta,\phi)
\nonumber
\end{eqnarray}
It is clear now that $\Pexpr(0,R)=\Czz$.

If  $f=\sum_{nm} f_{nm}Y^m_n\in\Pexpr(0,R)$ then $\pi_0(f)=f_{00}$.
However 
\begin{eqnarray}
\frac1{4\pi R^2}\int Y^m_n(\theta,\phi) \sin\theta d\phi d\theta &=&
\left\{ 
\begin{array}{ll}
1 & \qquad m=0 \hbox{ and } n=0 \cr
0 & \qquad \hbox{otherwise}
\end{array} \right.
\nonumber
\end{eqnarray}
so (\ref{lim_def_IP}).

\end{proof}


\section{The Moyal Bracket}
\label{ch_Moy}

The Moyal Bracket is defined by the limit of the commutator of two
elements
\begin{eqnarray}
&&\{\bullet,\bullet\}:\Czz(S^2)\times \Czz(S^2)
\mapsto
\Czz(S^2)
\cr
&&\{\fs f,\fs g\} =
\lim_{\kappa\to0}
\left(\frac{1}{i\kappa} 
[\Psi_{0,\kappa}(\fs f),\Psi_{0,\kappa}(\fs g)] \right)
\label{Moy_def_Moy}
\end{eqnarray}

\begin{theorem}
\label{Moy_thm}
If $\fs{f},\fs{g}\in \Czz(S^2)$ then
we have the the Moyal bracket is the natural bracket arising from the
symplectic form on $S^2$.
\begin{eqnarray}
\{\fs f,\fs g\} &=& 
\frac1{R\sin\theta}\left(
\frac{\partial \fs f}{\partial\phi} \, 
\frac{\partial \fs g}{\partial\theta} 
-
\frac{\partial \fs f}{\partial\theta} \, 
\frac{\partial \fs g}{\partial\phi} 
\right)
\label{Moy_Moy_res}
\end{eqnarray}
\end{theorem}


\begin{proof}{}
Since we are dealing with only finite sums of basis elements we need
not worry about limits. Since the Moyal bracket is linear in both
terms we need only consider its effect on basis elements.  Let
$f=\Psi_{0,\kappa}(\fs f)$ and $g=\Psi_{0,\kappa}(\fs g)$ be
eigenvectors of $e_z$.  For this proof we write
$\partial_\phi=\partial/\partial\phi$.  We note that for a polynomial
$p(\fs z)$
\begin{eqnarray}
\partial_\phi\Big(\fs{J}_+^m p(\fs{z})\Big) &=& -im\fs{J}_-^m
p(\fs{z})  
\cr
\partial_\phi\Big(\fs{J}_-^m p(\fs{z})\Big) &=& im\fs{J}_-^m p(\fs{z}) 
\cr
\frac1{R\sin\theta}\partial_\theta
\Big(\fs{J}_+^m p(\fs{z})\Big) &=&
-\fs{J}_+^m \left( 
p'(\fs{z}) -
\frac{m\fs{z}}{R^2-\fs{z}^2}p(\fs{z}) \right)
\cr
\frac1{R\sin\theta}\partial_\theta
\Big(\fs{J}_-^m p(\fs{z})\Big) &=&
-\fs{J}_-^m \left( 
p'(\fs{z}) -
\frac{m\fs{z}}{R^2-\fs{z}^2}p(\fs{z}) \right)
\nonumber
\end{eqnarray}
It is necessary to consider separately the cases that the eigenvalues
of $f$ and $g$ have (1) the same sign and (2) different signs. Let us
first consider the case when the eigenvalues of $f$ and $g$ have
positive sign. Then we can write
\begin{eqnarray}
f = J_+^a p(z) \hbox{ , and } g = J_+^b q(z) 
\nonumber
\end{eqnarray}
where $a,b\ge0$ and $p(z),q(z)$ are polynomials. Then
\begin{eqnarray}
[f,g] &=&
J_+^a p(z) J_+^b q(z) - J_+^b q(z) J_+^a p(z) \cr
&=&
J_+^{a+b} ( p(z+\kappa b) q(z) -  p(z) q(z+\kappa a) ) \cr
&=& 
J_+^{a+b} \Big( (p(z+\kappa b) - p(z))q(z) -
p(z)(q(z+\kappa a) -  q(z)) \Big) 
\nonumber
\end{eqnarray}
In the limit as $\kappa\to0$ we have $J_+\to\fs{J}_+$, 
$p(z)\to\fs{p}(z)$ and
$\nfrac1\kappa(p(z+b\kappa)-p(z))\to b\fs{p}'(z)$. So
\begin{eqnarray}
\lim_{\kappa\to0} \left(\frac{1}{i\kappa} [f,g] \right)
&=&
-i \fs{J}_+^{a+b}  
( b\fs{p}'(\fs z)\fs{q}(\fs z)-a\fs{p}(\fs z)\fs{q}'(\fs z) )\cr
\cr
&=&
(R\sin\theta)^{-1}
(\partial_\phi \fs f \, \partial_\theta \fs g
-
\partial_\theta \fs f \, \partial_\phi \fs g )
\nonumber
\end{eqnarray}
Hence true in this case.  If the eigenvalues of $f$ and $g$ both have
negative sign then we note that
\begin{eqnarray}
f^\dagger \to \cnj{\fs f}
\nonumber
\end{eqnarray}
and 
\begin{eqnarray}
[f,g]/(i\kappa) &=& 
[g^\dagger,f^\dagger]^\dagger/(i\kappa) 
\to
\cnj{\{\cnj{\fs f},\cnj{\fs g}\}} 
=
{\{{\fs f},{\fs g}\}} 
\nonumber
\end{eqnarray}
\vskip 1em

Now consider case (2) we write
\begin{eqnarray}
f = J_+^a p(z)\ ,\  
g = J_-^b q(z)
\hbox{ , and } 
\varrho_-^b(z) = J_+^bJ_-^b\ ,\ 
\varrho_+^b(z) = J_-^bJ_+^b
\nonumber
\end{eqnarray}
where $a,b\ge0$ and $p(z),q(z)$ are polynomials. Consider first $a\ge
b$
\begin{eqnarray}
[f,g] &=& 
J_+^a p(z) J_-^b q(z) - J_-^b q(z) J_+^a p(z) 
\cr
&=& \lgap
J_+^{a}J_-^{b} p(z-\kappa b) q(z) 
-  
J_-^{b}J_+^{a} p(z) q(z+\kappa a) ) 
\cr
&=& \lgap
J_+^{a-b}\varrho_-^b(z) p(z-\kappa b) q(z)
- 
J_+^{a-b}\varrho_+^b(z+(a-b)\kappa) p(z) q(z+\kappa a)  
\cr
&=&
J_+^{a-b} \bigg(
\begin{array}[t]{@{\lgap}l}
\varrho_-^b(z)(p(z-\kappa b)-p(z))q(z)
+ 
\varrho_-^b(z)p(z)(q(z)-q(z+\kappa a))
\cr
+
(\varrho_-^b(z) - \varrho_+^b(z))p(z) q(z+\kappa a)
+
(\varrho_+^b(z) - \varrho_+^b(z+\kappa(a-b))) p(z) q(z+\kappa a) \bigg)
\end{array}
\nonumber
\end{eqnarray}
In the limit $\kappa\to0$ this becomes
\begin{eqnarray}
\lefteqn{\lim_{\kappa\to0}
\left(\frac{1}{i\kappa} [f,g] \right)} 
\qquad\qquad
\cr
&=&
\fs J_+^{a-b} \bigg(
(R^2-\fs z^2)^b \fs p'(\fs z)(-b)\fs q(\fs z)
+
(R^2-\fs z^2)^b \fs p(\fs z)(-a)\fs q'(\fs z)
\cr &&
+
2(R^2-\fs z^2)^{b-1} b^2 \fs z \fs p(\fs z)\fs q(\fs z)
-
(a-b)\fs p(\fs z)\fs q(\fs z)\frac{d}{dz}(R^2-\fs z^2)^b \bigg)
\cr
&=&
\fs J_+^{a-b} \bigg(
(R^2-\fs z^2)^b ( -b \fs p'(\fs z)\fs q(\fs z) - 
a \fs p(\fs z)\fs q'(\fs z) )
+
(R^2-\fs z^2)^{b-1} 2zab \fs p(\fs z)\fs q(\fs z) \bigg)
\cr
&=&
(R\sin\theta)^{-1}\Big(
\partial_\phi(\fs J_+^a \fs p(\fs z))
\partial_\theta(\fs J_-^b \fs q(\fs z)) -
\partial_\phi(\fs J_-^b \fs q(\fs z))
\partial_\theta(\fs J_+^a \fs p(\fs z)) \Big)
\nonumber
\end{eqnarray}
Likewise if $b>a$ then we consider $[f^\dagger,g^\dagger]$ as before.
\end{proof} 


We now wish to consider how we can extend this theorem to cover the
largest possible subset of $L^2(S^2)$. For our case a sufficient
extension to $\Czz(S^2)$ is given by the set 
\begin{eqnarray}
\left\{ f=\sum_{nm}f_{nm}Y^m_n\in L^2(S^2) \ \bigg|\ 
|f_{nm}|\sim n^{-3} \right\}
\end{eqnarray}
This is because in this case $\partial_\phi f\in L^2(S^2)$ and
$(\sin\theta)^{-1}\partial_\theta f\in L^2(S^2)$.  Hence the right
hand side of (\ref{Moy_Moy_res}) is defined. 


\subsection*{Limit of the operators as $\kappa\to 0$}

As already mentioned the operations $e_z,e_\pm,\Delta$ from
(\ref{Pmn_def_ez}) to (\ref{Pmn_def_Del}) mean they identically vanish
if $\kappa=0$. We therefore calculate the first non-vanishing term in
there expansions. We see that $e_x,e_y,e_z,e_+,e_-$ are vector fields
whilst $\Delta$ is a second-order differential operator corresponding
to Laplacian.

\begin{theorem}
In the limit $\kappa\to0$ we have
\begin{eqnarray}
(i\kappa)^{-1} e_x &\to& 
\cos\phi\partial_\theta - \cot\theta \sin\phi \partial_\phi     
\\
(i\kappa)^{-1} e_y &\to& 
-\sin\phi\partial_\theta - \cot\theta \cos\phi \partial_\phi     
\\
(i\kappa)^{-1} e_z &\to& \partial_\phi 
\\
(i\kappa)^{-1} e_+ &\to& 
e^{-i\phi}\left( \partial_\theta -i\cot\theta\partial_\phi \right)
\\
(i\kappa)^{-1} e_- &\to&   
e^{i\phi}\left( \partial_\theta + i\cot\theta\partial_\phi \right)
\end{eqnarray}
and as one would expect
\begin{eqnarray}
(i\kappa)^{-1} e_x(\fs y) &\to& \fs z \textup{ and cyclic permutations
}
\end{eqnarray}
These are not independent since
\begin{eqnarray}
\fs x e_x+\fs y e_y+\fs z e_z=0
\end{eqnarray}
Also the ``fuzzy'' Laplace operator tends to the usual Laplace
operator on the sphere.
\begin{eqnarray}
-\kappa^{-2}\Delta &\to& 
\partial^2_\theta + \cot\theta\partial_\theta + (\sin\theta)^{-2}
\partial^2_\phi
\end{eqnarray}
\end{theorem}

\begin{proof}{}
The expressions for $e_+,e_-,e_z$ come by substituting $\fs J_+,\fs J_-,\fs z$ as
one of the term in the Moyal bracket. The other expressions are
derived from these.
\end{proof}

\section{Discussion}
\label{ch_disc}

This work forms a basis for the investigation into the differential
and connection structures on the fuzzy sphere. (Follow references in
\cite{Madore_book}). Since $\Pexpr^n$ is a $2n+1$ dimensional
representation of $su(2)$ we may consider these representing Bosonic
states. One should be able to create another basis of the Fermionic
states. This may look like $P^m_n$ with $m$ and $n$ positive odd
multiples of a half. 

It would be useful to know how these results can be extended for other
algebras. For $su(3)$ one would consider replacing $J_\pm$ with
$u_\pm$ and $v_\pm$, where $u_-,v_-$ are the root system. In this case
we would have a basis something like
\begin{eqnarray}
e_{u_-}^a e_{v_-}^b (u_+^c v_+^d) 
\nonumber
\end{eqnarray}
with some relation for the $a,b,c,d$. 
The case of $su(n)$ would be
equivalent using the root system. We might be able to extend this to
all Lie algebras of compact Lie group.

This article has demonstrated how one can use a quotient of a free
noncommuting algebra on a finite set of elements to examine a
geometry. For existence of an exterior algebra this quotient algebra
must form a ``generalised algebra'' \cite{Gratus4}.  This may be
necessary for quantising general manifolds.


\subsection*{Acknowledgements}

The author would like to thank John Madore, Luiz Saeger and Jihad
Mourad for useful discussions which motivated this work.  The author
would also like to thank the Royal Society of London for a European
Junior Fellowship, and Richard Kerner and the Laboratoire de
Gravitation et Cosmologie Relativistes, Paris~VI for their
hospitality.


\section*{Appendix A: Some Results about the Universal Enveloping
Algebra $\Sexpr$}

\setcounter{equation}{0}
\renewcommand{\theequation}{A\arabic{equation}}

This section is needed to establish the fact that
$P^m_n\in\Pexpr^n$. Since the dimension of $\Pexpr^n$ is $2n+1$ is is
obvious that $\{P^m_n,m=-n\ldots n\}$ form an orthonormal basis for
$\Pexpr^n$. (There may be an easier way without introducing this
machinery). Some of these results are mentioned without proof in
\cite{Bayen1}

Let us define the algebra $\Sexpr$ by
\begin{eqnarray}
\Sexpr &=&
\{ \mbox{Free noncommuting algebra of polynomials in $x,y,z$ } \}
\Big/\sim
\label{Sex_def_S}
\end{eqnarray}
where
\begin{eqnarray}
[x,y]\sim i\kappa z,\
[y,z]\sim i\kappa x,\
[z,x]\sim i\kappa y,\
\label{Sex_quot}
\end{eqnarray}
There is a natural basis of this algebra given by the 3-vector
$\Ssym(a,b,c)$ with $a,b,c\in\Intg$ and $a,b,c\ge 0$. This represent
the sum of all symmetric permutations of the word $x^ay^bz^c$ each with
coefficient $1$. Thus for example
\begin{eqnarray}
\Ssym(2,1,0) = x^2y + xyx + yx^2
\nonumber
\end{eqnarray}
It is easy to show that $\Ssym(a,b,c)$ has $(a+b+c)!/(a!b!c!)$ terms.
For consistency we define $\Ssym(a,b,c)=0$ if either $a<0$, $b<0$ or
$c<0$. These will be known as $\calS$-expressions. 

We can define the formal trace by 
\begin{eqnarray}
&& \ftr:\Sexpr\mapsto\Sexpr \cr\lgap
&& \ftr\left(\sum_{a_1\ldots a_p}
f_{a_1\ldots a_p}x_{a_1}\cdots x_{a_p}\right)
=
\sum_{b=1}^3\sum_{a_1\ldots a_p}
f_{b,b,a_3\ldots a_p}x_{a_3}\cdots x_{a_p}
\end{eqnarray}
Thus if $f\in\Sexpr$ and $\ftr(f)=0$ then $f\in\Pexpr$. We have the
following theorem for the manipulation of the $\calS$-expressions.

\begin{theorem}
The formal trace of an $\calS$-expression given by
\begin{eqnarray}
\tr\Big(\Ssym(a,b,c) \Big) &=&
\Ssym(a-2,b,c) + \Ssym(a,b-2,c) + \Ssym(a,b,c-2)
\label{Sex_form_tr}
\end{eqnarray}
The commutator of $x$ and an $\calS$-expression is given by
\begin{eqnarray}
e_x\Ssym(a,b,c) =
\left[x,\Ssym(a,b,c)\right]
=
-\kappa(b+1)\Ssym(a,b+1,c-1) + \kappa(c+1)\Ssym(a,b-1,c+1)
\label{Sex_adx}
\end{eqnarray}
and cyclic permutation for $e_y$, and $e_z$.
The relationship between these operations is given by
\begin{eqnarray}
\tr\circ e_x = e_x\circ\tr
\label{Sex_adx_tr}
\end{eqnarray}
and similarly for $e_y$, and $e_z$.
We can split an $\calS$-expression to the $m$ order to give
\begin{eqnarray}
\Ssym(a,b,c) &=&
\sum_{d+e+f=m} \Ssym(a-d,b-e,c-f)\Ssym(d,e,f)
\label{Sex_split_S}
\end{eqnarray}
\end{theorem}


\begin{proof}{}
We can think of the $\calS$-expression $\Ssym(a,b,c)$ as being a sum
of terms. Let $w$ be a permutation of $x^dy^ez^f$, with $d\le a,e\le
b,f\le c$.  Take all the terms in $\Ssym(a,b,c)$ which start with $w$.
These terms must finish with each term in $\Ssym(a-d,b-e,c-f)$. This
works with each permutation of $x^dy^ez^f$ so $\Ssym(a,b,c)$ must
contain the term $\Ssym(d,e,f)\Ssym(a-d,b-e,c-f)$. Now if we let
$d,e,f$ run over all sets $d+e+f=m$ and $d,e,f\ge0$ then this covers
all possibilities, and no two are repeated. Hence (\ref{Sex_split_S}).

From (\ref{Sex_split_S}) putting $m=2$ we have
\begin{eqnarray}
&&\Ssym(a,b,c) = 
\Ssym(2,0,0)\Ssym(a-2,b,c)+
\Ssym(0,2,0)\Ssym(a,b-2,c)+
\Ssym(0,0,2)\Ssym(a,b,c-2) 
\cr && \qquad +
\Ssym(1,1,0)\Ssym(a-1,b-1,c)+
\Ssym(1,0,1)\Ssym(a-1,b,c-1)+
\Ssym(0,1,1)\Ssym(a,b-1,c-1)
\nonumber
\end{eqnarray}
Now $\tr(x^2)=\tr(y^2)=\tr(z^2)=1$ hence (\ref{Sex_form_tr}).  Proof
of (\ref{Sex_adx}) is by induction on the order of the polynomial by
the use of (\ref{Sex_split_S}) with $m=1$. Proof is (\ref{Sex_adx_tr})
is by direct substitution.
\end{proof}

\begin{corol}
\label{Sex_corol}
We are now in a position to prove that $e_x:\Pexpr^n\mapsto\Pexpr^n$
and likewise for $e_y,e_z,e_\pm,\Delta$. Also $P^m_n\in\Pexpr^n$.
\end{corol}

\begin{proof}{of $P^m_n\in\Pexpr^n$}
From the definition of $P^n_n$
\begin{eqnarray}
P^n_n &=& \alpha_n J_+^n = \alpha_n(x+iy)^n \cr
&=&\alpha_n \sum_{r=0}^n i^r\Ssym(n-r,r,0)
\nonumber
\end{eqnarray}
also
\begin{eqnarray}
\ftr(P^n_n) &=& \alpha_n \sum_{r=0}^n i^r
(\Ssym(n-r-2,r,0)+\Ssym(n-r,r-2,0))
\cr
&=&
\alpha_n \sum_{r=0}^{n-2} i^r
\Ssym(n-r-2,r,0)
+
\alpha_n \sum_{r=0}^{n-2} i^{r+2}
\Ssym(n-r-2,r,0)
=0
\nonumber
\end{eqnarray}
So $P^n_n$ is a symmetric formally trace-free polynomial of order $n$.
So $P^n_n\in\Pexpr^n$. 

From (\ref{Sex_adx}) we see that if $f$ is an $\calS$-expression of
order $n$ then
so is $e_-(f)$. From (\ref{Sex_adx_tr}) we see that if $f\in\Pexpr^n$
then $e_-(f)\in\Pexpr^n$. So $P^m_n\in\Pexpr^n$.
\end{proof}




\begin{thebibliography}{99}

\bibitem{Berezin74} {Berezin, F.A.: Quantization, Math. USSR Izvestija
{\bf 8}, 1109--1165 (1974).}

\bibitem{Berezin75a}{Berezin, F.A.: General concept of quantization,
Commun. Math. Phys. {\bf 40}, 153--174 (1975).}
 
\bibitem{CaheGutt90}{Cahen, M., Gutt, S., Rawnsley, J.: Quantization
of K\"ahler Manifolds II, Transactions of the American Math. Soc. {\bf
337}, 73--98, (1993).}

\bibitem{Madore2} J. Madore, {\it The Commutative Limit of a Matrix
Geometry.} J. Math. Phys. {\bf 32}(2), 1991 332-335

\bibitem{Madore_book} J. Madore 1995, {\it An Introduction to
Noncommutative Differential Geometry and its Physical Applications},
Cambridge University Press.

\bibitem{klimcik1} H. Grosse, C. Klim\v c\'{\i}k, P. Pre\v{s}najder,
Simple Field Theoretical Models on Noncommutative Manifolds,
hep-th/9510177

\bibitem{Watamura1} U. Carow-Watamura, S. Watamura, Chirality and Dirac
Operators on the Noncommutative Sphere. Tohoku University Preprint 498
hep-th/9605003.

\bibitem{DeWit1} B. de Wit, U. Marquard, H. Nicolai. {\it
Area-Preserving Diffeomorphisms and Supermembrane Lorentz Invariance}
Commun. Math. Phys. 128, (1990) 39-62

\bibitem{Grosse1} H. Grosse, P. Pre\v{s}najder {\it The Construction of
Noncommutative Manifolds using Coherent states} Lett. Math. Phys. {\bf
28 }, 1993 239-250

\bibitem{Perelomov1} A. Perelomov 1986, {\it Generalized Coherent
States and their Applications} Springer-Verlag

\bibitem{Molin1} P. Molin, {\it Application Momentum Quantization}
Lett. Math. Phys. {\bf 25}, 1992 213-225

\bibitem{Molin2} P. Molin, {\it A star Product on the Spherical
Harmonics} Lett. Math. Phys. {\bf 38}, 1996 227-236

\bibitem{Bayen1} F. Bayen, C. Fronsdal, {\it Quantization on the
Sphere} J. Math. Phys. {\bf 22}(7), 1981 1345-1349

\bibitem{Nik} A. F. Nikiforov, S. K. Suslov, V. B. Uvarov, Classical
Orthogonal Polynomials of a Discrete Variable.

\bibitem{Gratus4} J. Gratus, {\it Non Commutative Differential
Geometry, and the matrix Representations of Generalised Algebras.}
to be published J. Geometry Phys. Accepted March 1997.

\end{thebibliography}
\end{document}